\documentclass{amsart}%
\usepackage{amsmath}
\usepackage{graphicx}
\usepackage{amsfonts}
\usepackage{amssymb}%
\newtheorem{theorem}{Theorem}
\theoremstyle{plain}

\newtheorem{definition}{Definition}
\newtheorem{example}{Example}

\newtheorem{lemma}{Lemma}

\newtheorem{proposition}{Proposition}
\newtheorem{remark}{Remark}

\numberwithin{equation}{section}
\begin{document}
\title{Calculating the Greeks by Cubature formulas}

\begin{abstract}
We provide cubature formulas for the calculation of derivatives of expected
values in the spririt of Terry Lyons and Nicolas Victoir. In financial
mathematics derivatives of option prices with respect to initial values, so
called Greeks, are of particular importance as hedging parameters. The proof
of existence of Cubature formulas for Greeks is based on an argument, which
leads to the calculation of Greeks in an asymptotic sense -- even without
H\"{o}rmander's condition. Cubature formulas then allow to calculate these
quantities very quickly. Simple examples are added to the theoretical exposition.

\end{abstract}
\author{Josef Teichmann}
\address{Technical University of Vienna, e105, Wiedner Hauptstrasse 8-10, A-1040 Wien, Austria}
\email{jteichma@fam.tuwien.ac.at}
\thanks{The author acknowledges the support from the RTN network HPRN-CT-2002-00281
and from the FWF grant Z-36.}
\subjclass[2000]{60H07, 35R60, 65C30}
\keywords{Cubature Formulas, Wiener Space, Malliavin Calculus, nilpotent Lie groups,
H\"ormander's Theorem}
\maketitle

\section{Introduction}

Cubature formulas provide approximative values for integrals with respect to a
given measure. The well-developed theory of cubature formulas in finite
dimensions was recently applied to provide Cubature formulas on Wiener space
by Terry Lyons and Nicolas Victoir in the beautiful seminal article
\cite{Lyo/Vic:04}.

We try to extend the framework of Terry Lyons and Nicolas Victoir to the
calculation of Greeks. We briefly outline in the introduction the main results
and possible applications of it. In the sequel we shall always work with
$C^{\infty}$-bounded vector fields $V_{i}$ and $C^{\infty}$-bounded function
$f$. Given a stochastic differential equation in $\mathbb{R}^{N}$ of the type%
\begin{align*}
dY_{t}^{y}  &  =V_{0}(V_{t}^{y})dt+\sum_{i=1}^{d}V_{i}(Y_{t}^{y})\circ
dB_{t}^{i},\\
Y_{0}^{y}  &  =y,
\end{align*}
on a stochastic basis $(\Omega,(\mathcal{F}_{t})_{0\leq t\leq T},P)$ with
$d$-dimensional Brownian motion, then cubature formulas provide a method to
approximate $E(f(Y_{t}^{y}))$, namely -- choosing a degree $m$ of accuracy --
there is number $r\geq1$ and there are $H^{1}$-trajectories $\omega
_{j}:[0,T]\rightarrow\mathbb{R}^{d+1}$ and weights $\lambda_{j}>0$,
$j=1,\dots,r$, such that%
\[
E(f(Y_{t}^{y}))=\sum_{j=1}^{r}\lambda_{j}f(Y_{t}^{y}(\omega_{j}))+\mathcal{O}%
(t^{\frac{m+1}{2}}),
\]
where $Y_{t}^{y}(\omega_{j})$ denotes the point in $\mathbb{R}^{N}$, which is
obtained by solving the ordinary, non-autonomous differential equation%
\begin{align*}
\frac{dY_{t}^{y}(\omega_{j})}{dt}  &  =\sum_{i=0}^{d}V_{i}(Y_{t}^{y}%
(\omega_{j}))\frac{d\omega_{i}^{i}}{dt},\\
Y_{0}^{y}(\omega_{j})  &  =y.
\end{align*}
The constant in the order estimate depends in general on derivatives up to
order $m+1$, in the hypo-ellitptic case this dependence can be improved to the
first derivative. Note also the important fact that the number $r$ of
trajectories can be bounded by a number depending on $d$ and $m$, but not on
$N$ (see \cite{Lyo/Vic:04} for the precise estimate). Obviously this procedure
can be iterated with small time steps by the Markov property, which yields
then a high order numerical approximation scheme for the calculation of
expected values $E(f(Y_{t}^{y}))$.

In the hypo-elliptic case we can calculate the derivatives of $y\mapsto
E(f(Y_{t}^{y}))$ with respect to the initial value by weight formulas, i.e.
for any direction $v\in\mathbb{R}^{N}$ there is a random variable $\pi$ such
that%
\begin{equation}
\frac{d}{d\epsilon}|_{\epsilon=0}E(f(Y_{t}^{y+\epsilon v})))=E(f(Y_{t}^{y}%
)\pi) \label{malliavin formula}%
\end{equation}
for bounded measurable $f$. Explicit formulas for $\pi$ can be given by
Malliavin Calculus (see for instance \cite{GobMun:02} and \cite{TeiTou:04})
and are in principle well-known since long time (see for instance
\cite{Bis:84}).

If we fix an order of approximation $m\geq1$ and a (homogenous) direction $v$,
such that (for the notations see Section \ref{cubature})%
\[
v=\sum_{\substack{I\in\mathcal{A}\setminus\{\emptyset,(0)\}\\\deg(I)=k}%
}v_{I}[\cdots\lbrack V_{i_{1},}V_{i_{2}}],V_{i_{3}}],\dots,V_{i_{k}}](y)
\]
holds with some $1\leq k\leq m$, then we can construct a universal weight
$\pi_{d,1}^{m}$ (calculated in a \emph{universal way} from $v_{I}$ and the
Brownian motions, see Remark \ref{universal remark}), such that%
\[
\frac{d}{d\epsilon}|_{\epsilon=0}E(f(Y_{t}^{y+\epsilon v})))t^{\frac{k}{2}%
}=E(f(Y_{t}^{y})\pi_{d,1}^{m})+\mathcal{O}(t^{\frac{m+1}{2}})
\]
holds true, where the constant in the order estimate only depends on the first
derivative of $f$. Theorem \ref{cubature for greeks} asserts that we are able
to find weights $\mu_{j}\neq0$ and $H^{1}$-trajectories $\omega_{j}%
:[0,T]\rightarrow\mathbb{R}^{d+1}$ such that
\begin{equation}
\frac{d}{d\epsilon}|_{\epsilon=0}E(f(Y_{t}^{y+\epsilon v}))t^{\frac{k}{2}%
}=\sum_{j=1}^{r}\mu_{j}f(Y_{t}^{y}(\omega_{j}))+\mathcal{O}(t^{\frac{m+1}{2}})
\label{formula for greeks}%
\end{equation}
holds true. Also here the constant in the order estimate depends on
derivatives of $f$ up to order $m+1$. The number $r$ can be estimated
similarly as in \cite{Lyo/Vic:04}.

A combination of the original Cubature formula and formula
(\ref{formula for greeks}) then yields by iteration a high order numerical
approximation scheme, i.e.%
\begin{align}
\frac{d}{d\epsilon}|_{\epsilon=0}E(f(Y_{t}^{y+\epsilon v}))  &  =\frac
{d}{d\epsilon}|_{\epsilon=0}E(E(f(Y_{t}^{y+\epsilon v})|\mathcal{F}_{s_{0}%
}))\label{recipe1}\\
&  =\frac{d}{d\epsilon}|_{\epsilon=0}E(E(f(Y_{t-s_{0}}^{z}))|_{z=Y_{s_{0}%
}^{y+\epsilon v}}). \label{recipe2}%
\end{align}
Calculation of the inner expectation by the original cubature formula and of
the outer expectation by formula (\ref{formula for greeks}) lead to the
desired iterative procedure: in the uniformly hypo-elliptic case (see
\cite{Lyo/Vic:04} for details) we consequently obtain the following recipe.

\begin{itemize}
\item choose a small $s_{0}>0$ and relatively large $t-s_{0}$. Subdivide
$s=t_{0}<t_{1}<\dots<t_{k}=t$ with $s_{i}:=t_{i}-t_{i-1}$ for $i=1,\dots,k$.

\item calculate by (\ref{formula for greeks}) with the function $z\mapsto
E(f(Y_{t-s_{0}}^{z}))$ for a certain degree of accuracy $m$. The order
estimate yields (by careful choice of homogenous) directions $\mathcal{O}%
(s^{\frac{m+1-k}{2}})$ for some $1\leq k\leq m$. The constant in the order
estimate depends on the $(m+1)$st derivative of $z\mapsto E(f(Y_{t-s}^{z}))$,
which -- by methods of Malliavin Calculus -- can be reduced to dependence on
the first derivative of $f$, i.e. there is an absolute bound given by%
\[
K||f||_{Lip}s_{0}^{\frac{m+1-k}{2}},
\]
since $t-s_{0}$ is bounded away from $0$.

\item calculate $z\mapsto E(f(Y_{t-s}^{z}))$ by Cubature formulas with certain
degree of accuracy $m^{\prime}$ (see \cite{Lyo/Vic:04}). The order estimate,
which again -- by methods from Malliavin Calculus (see \cite{Kus:01}) -- can
be reduced to dependence on the first derivative of $f$, i.e. there is an
absolute bound given by
\[
K||f||_{Lip}\left(  \sum_{i=1}^{k-1}\frac{s_{i}^{\frac{m^{\prime}+1}{2}}%
}{(t-t_{i})^{\frac{m^{\prime}}{2}}}+s_{k}^{\frac{1}{2}}\right)  .
\]

\item finally we obtain a high order approximation scheme for $\frac
{d}{d\epsilon}|_{\epsilon=0}E(f(Y_{t}^{y+\epsilon v}))$ through solving
ordinary differential equations with estimate of the absolute error through%
\[
K||f||_{Lip}\left(  s_{0}^{\frac{m+1-k}{2}}+\sum_{i=1}^{k-1}\frac{s_{i}%
^{\frac{m^{\prime}+1}{2}}}{(t-t_{i})^{\frac{m^{\prime}}{2}}}+s_{k}^{\frac
{1}{2}}\right)  .
\]
The procedure only involves ordinary differential equations and knowledge on
the cubature trajectories. By choosing $m$ big enough, one can bound the
difference $m-k$ from below.
\end{itemize}

\begin{remark}
A direct application of formula (\ref{malliavin formula}) is often very slow
numerically in the purely hypo-elliptic case (see \cite{GobMun:02}). In the
elliptic case there are very efficient ways to simulate $\pi$ (see
\cite{FouLasLebLioTou:99}). Hence in particular the purely hypo-elliptic case
constitutes an interesting field of applications, since one does not need to
worry about invertibility of the covariance matrix.
\end{remark}

\begin{remark}
For all cubature formulas there are algebraic methods to determine the
trajectories $\omega_{j}$ and the weights $\lambda_{j}$, or $\mu_{j}$,
respectively. These methods are based on calculus in nilpotent Lie groups,
which constitute \emph{the} local model for hypo-elliptic diffusions in the
sense of Gromov (see for instance the very well written monograph
\cite{Mon:02} for a general outline, and \cite{Bau:04}, \cite{Lyo:98} for
applications in the theory of -- stochastic -- differential equations).
\end{remark}

\section{Cubature on Wiener Space\label{cubature}}

Cubature formulas on Wiener space (see the seminal article \cite{Lyo/Vic:04})
rely on the analysis of hypo-elliptic diffusions on free, nilpotent groups on
the one hand (see the recent book \cite{Bau:04}), on the other hand on
stochastic Taylor expansion (see for instance \cite{Ben:89}). Furthermore the
proof of a Cubature formula on Wiener space appears as an application of
Tchakaloff's Theorem on finite dimensional cubature formulas (see
\cite{Put:97}). We state these three main ingredients and show -- in order to
explain the mathematical background of Cubature formulas -- a sketch of the
proof (following \cite{Lyo/Vic:04}).

We fix a probability space $(\Omega,(\mathcal{F}_{t})_{0\leq t\leq T},P)$
together with a $d$-dimensional Brownian motion $(B_{s}^{i})_{0\leq t\leq
T,i=1,\dots,d}$. We need the following notations for convenience:

\begin{enumerate}
\item We abbreviate by $\mathcal{A}$ the set of all finite sequences
$I:=(i_{1},\dots,i_{k})\in\{0,\dots,d\}^{k}$, $k\geq0$ and we define a degree
function on $\mathcal{A}$ by%
\[
\deg(i_{1},\dots,i_{k}):=k+\#\{l\text{ with }i_{l}=0\},
\]
which simply means that appearing $0$s are counted twice. Additionally, we
define $\deg(\emptyset)=0$. We define a semi-group structure on $\mathcal{A}$
via%
\begin{align*}
(i_{1},\dots,i_{k})\ast(j_{1},\dots,j_{l})  &  :=(i_{1},\dots,i_{k}%
,j_{1},\dots,j_{l}),\\
\emptyset\ast(i_{1},\dots,i_{k})  &  =(i_{1},\dots,i_{k})\ast\emptyset
=(i_{1},\dots,i_{k}).
\end{align*}

\item We denote by $\mathbb{A}_{d,1}^{m}$the free, nilpotent algebra with
$d+1$ generators $e_{0},\dots,e_{d}$, i.e. the set of all non-commutative
polynomials in those variables, such that the following nilpotency relations
hold: if $\deg(I)>m$, then $e_{i_{1}}e_{i_{2}}\cdots e_{i_{k}}=0$ for all
$I\in\mathcal{A}$. Hence $\mathbb{A}_{d,1}^{m}$ is a finite dimensional,
non-commutative, real algebra with unit element $1$ and we are given a grading
via the degree functions, i.e. a monomial $e_{i_{1}}e_{i_{2}}\cdots e_{i_{k}}$
is said to have degree $n$ if $\deg(i_{1},\dots,i_{k})=n.$ Denote by $W_{n}$
the linear span of all monomials of degree $n$, then we obtain%
\[
\mathbb{A}_{d,1}^{m}=\oplus_{n=0}^{m}W_{n},
\]
furthermore $W_{p}W_{q}\subset W_{p+q}$ (where we define $W_{p}=0$ for $p>m$
due to the given relations), $W_{0}=\mathbb{R}\cdot1$, so $\mathbb{A}%
_{d,1}^{m}$ is a graded algebra. We denote the canonical projections of $x$ on
the subspaces $W_{n}$ by $x_{n}$ for $n\geq0$.

\item On the finite dimensional algebra $\mathbb{A}_{d,1}^{m}$ we define the
exponential series%
\[
\exp(x):=\sum_{i=0}^{\infty}\frac{x^{i}}{i!},
\]
where the series converges everywhere due to the nilpotency relations. We
define the logarithm on elements $x$ with $x_{0}\neq0$. We identify $x_{0}$
with a real number, so the series%
\[
\log(x)=\log(x_{0})+\sum_{i=1}^{\infty}\frac{(-1)^{i-1}}{i}(\frac{x-x_{0}%
}{x_{0}})^{i}.
\]
is well-defined, since $\frac{x-x_{0}}{x_{0}}\in\oplus_{n=1}^{m}W_{n}$ (i.e.
the series is finite).

\item We define the Lie algebra $\mathfrak{g}_{d,1}^{m}$ generated by
$e_{0},\dots,e_{d}$ with respect to the Lie bracket $[x,y]:=xy-yx$. The Lie
algebra inherits the grading from the algebra $\mathbb{A}_{d,1}^{m}$ via
$U_{n}:=\mathfrak{g}_{d,1}^{m}\cap W_{n}$, for $n\geq0$. Hence%
\[
\mathfrak{g}_{d,1}^{m}=\oplus_{n=1}^{m}U_{n}%
\]
is a graded Lie algebra. In fact the Lie algebra is free, nilpotent of step
$m$, with $d$ generators of degree $1$ and one generator of degree $2$.

\item We denote the exponential image of $\mathfrak{g}_{d,1}^{m}$ by
$G_{d,1}^{m}:=\exp(\mathfrak{g}_{d,1}^{m})\subset\mathbb{A}_{d,1}^{m}$ and
call it the free, nilpotent Lie group. $G_{d,1}^{m}$ is indeed a Lie group as
closed subgroup of Lie group $1+\oplus_{n=1}^{m}W_{n}$. The tangent space at
$x\in G_{d}^{m}$ is spanned by the left (or right) translations $xw$ for
$w\in\mathfrak{g}_{d,1}^{m}$.

\item We need the canonical dilatations on Lie algebra and Lie group. We
define for $x\in\mathbb{A}_{d,1}^{m}$ an algebra homomorphism $\Delta_{t}$ via%
\[
\Delta_{t}(x_{n})=t^{n}x_{n}%
\]
for $x_{n}\in W_{n}$, $n\geq1$ and $t>0$. The homomorphism is well defined and
restricts to a Lie algebra homomorphism $\Delta_{t}$ on $\mathfrak{g}%
_{d,1}^{m}$ and and a Lie group homomorphism $\Delta_{t}$ for $t>0$.

\item We apply the notation $\circ dB_{t}^{0}:=dt$. We abbreviate the left
invariant vector field associated to $e_{i}$ by $D_{i}$, i.e. $D_{i}%
(x)=xe_{i}$ for $x\in G_{d,1}^{m}$. We define a stochastic process $(X_{t}%
^{x})_{0\leq t\leq T}$ on $G_{d,1}^{m}$ via%
\begin{align}
dX_{t}^{x}  &  =\sum_{i=0}^{d}D_{i}(X_{t}^{x})\circ dB_{t}^{i}%
\label{sde-liegroup}\\
&  =\sum_{i=0}^{d}X_{t}^{x}e_{i}\circ dB_{t}^{i},\\
X_{0}^{x}  &  =x,
\end{align}
and see immediately that $X_{t}^{x}=xX_{t}^{1}$ for $0\leq t\leq T$ almost surely.

\item Note here and in the sequel that vector fields on a vector space
$V:\mathbb{R}^{N}\rightarrow\mathbb{R}^{N}$ are used with double meaning:
either as tangent directions on a smooth geometric object, or as first order
differential operators on smooth functions $f:\mathbb{R}^{N}\rightarrow
\mathbb{R,}$%
\[
(Vf)(y):=df(x)\cdot V(x)
\]
for $x\in\mathbb{R}^{N}$. For smooth maps $G:\mathbb{R}^{N}\rightarrow
\mathbb{R}^{M}$ we apply the notion of the \emph{tangent map or Jacobian}%
\begin{equation}
dG(x)\cdot v:=\frac{d}{d\epsilon}|_{\epsilon=0}G(x+\epsilon v)
\label{jacobian}%
\end{equation}
for $x,v\in\mathbb{R}^{N}$. A vector field is called \emph{$C^{\infty}%
$-bounded} if all derivatives of order greater of order greater than $0 $ are bounded.
\end{enumerate}

In order to see the relation between free, nilpotent Lie groups on the one
hand and asymptotic analysis as $t\downarrow0$ on the other hand, we formulate
stochastic Taylor expansion (see \cite{Bau:04} and \cite{Ben:89}), which
simply results from iterating the defining equations for a stochastic
differential equation.

\begin{theorem}
[stochastic Taylor expansion]Given $C^{\infty}$-bounded vector fields
\[
V_{0},\dots,V_{d} : \mathbb{R}^{N}\rightarrow\mathbb{R}^{N}%
\]
on $\mathbb{R}^{N}$, then the diffusion process%
\begin{align*}
dY_{t}^{y}  &  =\sum_{i=0}^{d}V_{i}(Y_{t}^{y})\circ dB_{t}^{i},\\
Y_{0}^{y}  &  =y,
\end{align*}
for $y\in\mathbb{R}^{N}$ admits the following series expansion: for $n\geq0$
and $f\in C_{b}^{\infty}(\mathbb{R}^{N})$ the expansion%
\begin{align*}
f(Y_{t}^{y})  &  =\sum_{\substack{I\in\mathcal{A}\\\deg(I)\leq n}}V_{i_{1}%
}\cdots V_{i_{k}}f(y)\int_{0\leq t_{1}\leq\dots\leq t_{k}\leq t}\circ
dB_{t_{1}}^{i_{1}}\circ\cdots\circ dB_{t_{k}}^{i_{k}}+\\
&  +R_{n}(t,f),\\
R_{n}(t,f)  &  =\sum_{\substack{(i_{1},\dots,i_{k})\in\mathcal{A}\\\deg
(i_{1},\dots,i_{k})\leq n\\(i_{0},\dots,i_{k})\in\mathcal{A}\\\deg(i_{0}%
,\dots,i_{k})>n}}\int_{0\leq t_{0}\leq\dots\leq t_{k}\leq t}V_{i_{0}}\cdots
V_{i_{k}}f(Y_{t_{0}}^{y})\circ dB_{t_{0}}^{i_{0}}\circ\cdots\circ dB_{t_{k}%
}^{i_{k}},
\end{align*}
holds true for $0\leq t\leq T$. Additionally $\sqrt{E(R_{n}(t,f)^{2}%
)}=\mathcal{O}(t^{\frac{m+1}{2}})$, where the constant in the order estimate
depends on the $(m+1)$st derivative of $f$.
\end{theorem}

\begin{remark}
The above procedure works for ''any''\ differential equation (see the
literature on rough paths, for instance the seminal work \cite{Lyo:98}), in
particular, if we consider the non-autonomous, ordinary differential equation%
\[
dZ_{t}^{y}=\sum_{i=0}^{d}V_{i}(Z_{t}^{y})d\omega^{i}(t)
\]
for a $H^{1}$-trajectory $\omega:[0,T]\rightarrow\mathbb{R}^{d+1}$, we have
for $n\geq0$ and $f\in C_{b}^{\infty}(\mathbb{R}^{N})$ the expansion%
\begin{align*}
f(Z_{t}^{y})  &  =\sum_{\substack{I\in\mathcal{A}\\\deg(I)\leq n}}V_{i_{1}%
}\cdots V_{i_{k}}f(y)\int_{0\leq t_{1}\leq\dots\leq t_{k}\leq t}d\omega
^{i_{1}}(t_{1})\cdots d\omega^{i_{k}}(t_{k})+\\
&  +R_{n}(t,f),\\
R_{n}(t,f)  &  =\sum_{\substack{(i_{1},\dots,i_{k})\in\mathcal{A}\\\deg
(i_{1},\dots,i_{k})\leq n\\(i_{0},\dots,i_{k})\in\mathcal{A}\\\deg(i_{0}%
,\dots,i_{k})>n}}\int_{0\leq t_{0}\leq\dots\leq t_{k}\leq t}V_{i_{0}}\cdots
V_{i_{k}}f(Z_{t_{0}}^{y})d\omega^{i_{0}}(t_{0})\cdots d\omega^{i_{k}}(t_{k}),
\end{align*}
holds true for $0\leq t\leq T$. Notice here that in order to obtain an
asymptotic expansion (for a given fixed trajectory $\omega$) of certain order
in time $t$, one has to change the degree function.\label{non-auto taylor}
\end{remark}

\begin{example}
If we apply this series expansion to the process $X_{t}^{x}$ in the vector
space $\mathbb{A}_{d,1}^{m}$, we obtain for any linear function $\lambda
:\mathbb{A}_{d,1}^{m}\rightarrow\mathbb{R}$%
\[
\lambda(X_{t}^{x})=\sum_{\substack{I\in\mathcal{A}\\\deg(I)\leq m}%
}\lambda(xe_{i_{1}}\cdots e_{i_{k}})\int_{0\leq t_{1}\leq\dots\leq t_{k}\leq
t}\circ dB_{t_{1}}^{i_{1}}\circ\cdots\circ dB_{t_{k}}^{i_{k}}%
\]
for $x\in\mathbb{A}_{d,1}^{m}$ and $0\leq t\leq T$, since%
\[
D_{i}\lambda(x):=\frac{d}{d\epsilon}|_{\epsilon=0}\lambda(x+\epsilon
xe_{i})=\lambda(xe_{i})
\]
for $e_{i}\in\mathfrak{g}_{d,1}^{m}$, hence%
\[
(D_{i_{1}}\cdots D_{i_{k}}\lambda)(x)=\lambda(xe_{i_{1}}\cdots e_{i_{k}}).
\]
Consequently we obtain in $\mathbb{A}_{d,1}^{m}$%
\[
X_{t}^{1}=\sum_{\substack{I\in\mathcal{A}\\\deg(I)\leq m}}e_{i_{1}}\cdots
e_{i_{k}}\int_{0\leq t_{1}\leq\dots\leq t_{k}\leq t}\circ dB_{t_{1}}^{i_{1}%
}\circ\cdots\circ dB_{t_{k}}^{i_{k}}%
\]
for $0\leq t\leq T$, almost surely, by duality. This formula provides a nice
stochastic representation of the solution of equation \ref{sde-liegroup},
where we apply that $G_{d,1}^{m}$ is embedded in the algebra $\mathbb{A}%
_{d,1}^{m}$.
\end{example}

Stochastic differential equations provide solutions for the associated heat
equation via expectations, so we obtain for $f\in C_{b}^{\infty}%
(\mathbb{A}_{d,1}^{m})$,%
\[
(D_{0}+\frac{1}{2}\sum_{i=1}^{d}D_{i}^{2})(x\mapsto E(f(X_{t}^{x}%
)))=\frac{\partial}{\partial t}E(f(X_{t}^{x})),
\]
for all $x\in\mathbb{A}_{d,1}^{m}$. The heat equation itself admits a similar
construction as stochastic Taylor expansion, namely any solution $u(t,x)$ with
initial value $u(0,x)=f(x)$, where $f$ is $C^{\infty}$-bounded, can be written
as%
\begin{align*}
u(t,x)  &  =\sum_{j=0}^{n}\frac{t^{j}}{j!}(D_{0}+\frac{1}{2}\sum_{i=1}%
^{d}D_{i}^{2})^{j}f(x)+R_{n}(t,f),\\
R_{n}(t,f)  &  =\int_{0}^{t}\cdots\int_{0}^{s_{1}}(D_{0}+\frac{1}{2}\sum
_{i=1}^{d}D_{i}^{2})^{n+1}u(s_{0},x)ds_{0}\dots ds_{n}.
\end{align*}
In the case of $f=\lambda$ for a linear function $\lambda:\mathbb{A}_{d,1}%
^{m}\rightarrow\mathbb{R}$, this leads for the process $X_{t}^{x}$ in
$\mathbb{A}_{d,1}^{m}$, to the nice formulas%
\begin{align*}
E(\lambda(X_{t}^{x}))  &  =\sum_{\substack{I\in\mathcal{A}\\\deg(I)\leq
m}}\lambda(xe_{i_{1}}\cdots e_{i_{k}})E(\int_{0\leq t_{1}\leq\dots\leq
t_{k}\leq t}\circ dB_{t_{1}}^{i_{1}}\circ\cdots\circ dB_{t_{k}}^{i_{k}}),\\
u(t,x)  &  =\lambda(x\exp(t(e_{0}+\frac{1}{2}\sum_{i=1}^{d}e_{i}^{2}))),
\end{align*}
and hence%
\begin{equation}
\exp(t(e_{0}+\frac{1}{2}\sum_{i=1}^{d}e_{i}^{2}))=\sum_{\substack{I\in
\mathcal{A}\\\deg(I)\leq m}}e_{i_{1}}\cdots e_{i_{k}}E(\int_{0\leq t_{1}%
\leq\dots\leq t_{k}\leq t}\circ dB_{t_{1}}^{i_{1}}\circ\cdots\circ dB_{t_{k}%
}^{i_{k}}), \label{expectationformula}%
\end{equation}
for $0\leq t\leq T$, which is one basic formula for the construction of
Cubature formulas (see \cite{Bau:04} and \cite{Lyo/Vic:04}).

In order to obtain a cubature formula we need a version of Tchakaloff's
theorem (see \cite{Put:97} for details in the case where the support of $\mu$
is not compact):

\begin{theorem}
Let $\mu$ be a measure on $\mathbb{R}^{N}$ such that moments of all orders
exist. Fix a number $m$, then there are points $x_{1},\dots,x_{r}%
\in\operatorname*{supp}(\mu)$ and weights $\lambda_{j}>0$ for $j=1,\dots,r$,
such that%
\[
\int_{\mathbb{R}^{N}}f(y)\mu(dy)=\sum_{j=1}^{r}f(x_{j})\lambda_{j}%
\]
for all polynomials on $\mathbb{R}^{N}$ up to order $m$. Furthermore
$r\leq\dim_{\mathbb{R}}\operatorname*{Pol}_{m}(\mathbb{R}^{N})$, where
$\operatorname*{Pol}_{m}(\mathbb{R}^{N})$ is the vector space of polynomials
on $\mathbb{R}^{N}$ with degree less or equal $m$.
\end{theorem}

On basis of these prerequisites we can define cubature formulas on Wiener space:

\begin{definition}
Fix $m\geq1$ and $0<t\leq T$, a cubature formula on Wiener space is given by a
finite number of points $x_{1},\dots,x_{r}\in G_{d,1}^{m}$ and finitely many
weights $\lambda_{1},\dots,\lambda_{r}>0$, such that%
\[
E(X_{t}^{1})=\sum_{j=1}^{r}\lambda_{j}x_{j},
\]
or equivalently due to formula (\ref{expectationformula})%
\[
\exp(t(e_{0}+\frac{1}{2}\sum_{i=1}^{d}e_{i}^{2}))=\sum_{j=1}^{r}\lambda
_{j}x_{j}.
\]
Furthermore $r\leq\dim_{\mathbb{R}}\mathbb{A}_{d,1}^{m}$.
\end{definition}

Applying Tchakaloff's theorem to the law of the process $X_{t}^{1}$ in
$\mathbb{A}_{d,1}^{m}$, $0<t\leq T$, which is supported in $G_{d,1}^{m}$,
yields that we find points $x_{1},\dots,x_{r}\in\operatorname*{supp}%
((X_{t}^{1})_{\ast}P)\subset G_{d,1}^{m}$ and weights $\lambda_{j}>0$, such
that%
\[
E(X_{t}^{1})=\sum_{\substack{I\in\mathcal{A}\\\deg(I)\leq m}}e_{i_{1}}\cdots
e_{i_{k}}E(\int_{0\leq t_{1}\leq\dots\leq t_{k}\leq t}\circ dB_{t_{1}}^{i_{1}%
}\circ\cdots\circ dB_{t_{k}}^{i_{k}})=\sum_{j=1}^{r}\lambda_{j}x_{j},
\]
whence the existence of Cubature formulas for any $m\geq1$.

Finally we apply Chow's theorem of sub-riemannian geometry (see \cite{Bau:04}
and \cite{Mon:02}), which tells that every point $x\in G_{d,1}^{m}$ can be
reached by a $H^{1}$-horizontal curve, i.e. for every $x\in G_{d,1}^{m}$ we
find a $H^{1}$-curve $\omega:[0,T]\rightarrow\mathbb{R}^{d+1}$, such that the
solution of the non-autonomous ordinary differential equation%
\[
dZ_{t}^{1}(\omega)=\sum_{i=0}^{d}D_{i}(Z_{t}^{1}(\omega))d\omega^{i}(t)
\]
reaches $x$ at time $t>0$, i.e. $Z_{t}^{1}(\omega)=x$.

Taking for each of the $x_{j}$ appropriate curves $\omega_{j}$ with the
previous property, we obtain -- by Taylor expansion of the solution of the
non-autonomuous equation as in Remark \ref{non-auto taylor} -- the familiar
version of cubature formulas on Wiener space, namely%
\begin{gather*}
\sum_{\substack{I\in\mathcal{A}\\\deg(I)\leq m}}e_{i_{1}}\cdots e_{i_{k}%
}E(\int_{0\leq t_{1}\leq\dots\leq t_{k}\leq t}\circ dB_{t_{1}}^{i_{1}}%
\circ\cdots\circ dB_{t_{k}}^{i_{k}})=\\
\sum_{j=1}^{r}\lambda_{j}\sum_{\substack{I\in\mathcal{A}\\\deg(I)\leq
m}}e_{i_{1}}\cdots e_{i_{k}}\int_{0\leq t_{1}\leq\dots\leq t_{k}\leq t}%
d\omega_{j}^{i_{1}}(t_{1})\cdots d\omega_{j}^{i_{k}}(t_{k}).
\end{gather*}
This leads to a redefinition of cubature formulas, where we replace the points
$x_{i}$ by endpoints of evolutions of ordinary differential equations
$Z_{t}^{1}(\omega_{i})=x_{i}$ for $i=1,\dots,r$.

\begin{remark}
As proposed in \cite{Lyo/Vic:04} we could also find for any point $x\in
G_{d,1}^{m}$ a $H^{1}$-trajectory $\widetilde{\omega}:[0,T]\rightarrow
\mathbb{R}^{d}$, such that $\omega(t):=(t,\widetilde{\omega}(t))$ for
$t\in\lbrack0,T]$ satisfies $Z_{t}^{1}(\omega)=x$, hence we could equally work
with this simpler class of curves, which we do not do in this article.
\end{remark}

\begin{theorem}
\label{cubature on wiener space}Given $C^{\infty}$-bounded vector fields
$V_{0},\dots,V_{d}$ on $\mathbb{R}^{N}$, then the diffusion process%
\begin{align*}
dY_{t}^{y}  &  =\sum_{i=0}^{d}V_{i}(Y_{t}^{y})\circ dB_{t}^{i},\\
Y_{0}^{y}  &  =y,
\end{align*}
for $y\in\mathbb{R}^{N}$ admits the following cubature formula of degree
$m\geq1$: for $0<t\leq T$ we find $H^{1}$-curves $\omega_{1},\dots,\omega
_{r}:[0,T]\rightarrow\mathbb{R}^{d+1}$ and weights $\lambda_{1},\dots
,\lambda_{r}>0$, such that%
\[
E(f(Y_{t}^{y}))=\sum_{j=1}^{r}\lambda_{j}f(Y_{t}^{y}(\omega_{j}))+\mathcal{O}%
(t^{\frac{m+1}{2}}),
\]
for $f\in C_{b}^{\infty}(\mathbb{R}^{N})$. The constant for the order estimate
depends in general on derivatives of $f$ up to order $m+1$. The curve
$(Y_{t}^{y}(\omega))_{0\leq t\leq T}$ is understood as solution of%
\begin{align*}
dY_{t}^{y}(\omega)  &  =\sum_{i=0}^{d}V_{i}(Y_{t}^{y}(\omega))d\omega
^{i}(t),\\
Y_{0}^{y}(\omega)  &  =y.
\end{align*}

\end{theorem}

\begin{remark}
Notice that in the hypo-elliptic case we can improve Theorem
\ref{cubature on wiener space} by methods proved in \cite{Kus:01}. Fix
$m\geq1$. A family of random variables $\{Z_{I}$ for $I\in\mathcal{A}%
,\deg(I)\leq m\}$ is called $m$-moment generating if $Z_{\emptyset}=1$ and%
\[
E(Z_{I_{1}}Z_{I_{2}}\cdots Z_{I_{l}})=E(B^{\circ I_{1}}(1)\cdots B^{\circ
I_{l}}(1))
\]
for all $I_{j}\in\mathcal{A}$ with $\deg(I_{1})+\dots+\deg(I_{l})\leq m$. Here
we apply%
\begin{align*}
B^{\circ\emptyset}(t)  &  =1\\
B^{\circ I\ast(i)}(t)  &  :=\int_{0}^{t}B^{\circ I}(s)\circ dB_{s}^{i}%
\end{align*}
for $0\leq t\leq T.$ The results on order estimates in \cite{Kus:01} apply for
$m$-moment-generating families.

Given now a cubature formula, the coefficients $\int_{0\leq t_{1}\leq\dots\leq
t_{k}\leq1}d\omega_{j}^{i_{1}}(t_{1})\cdots d\omega_{j}^{i_{k}}(t_{k})$ in the
basis $e_{i_{1}}\cdots e_{i_{k}}$, which appear with probability $\lambda_{j}$
constitute a $m$-moment generating family (see \cite{Lyo/Vic:04}): this is
basically due to the fact that products of iterated integrals like $B^{\circ
I_{1}}(1)\cdots B^{\circ I_{l}}(1)$ can be decomposed as linear combination of
iterated integrals (for instance $B_{1}^{1}B_{1}^{2}=B^{\circ(1,2)}%
(1)+B^{\circ(2,1)}(1)$), hence the assertion on expectations of products
follows from the hypothesis on the expectation of iterated integrals.
Equivalently one can see this phenomenon also by transforming the defining
differential equation for the process $(X_{t}^{x})_{t\geq0}$ to the Lie
algebra in the exponential chart. It is remarkable that -- due to this
observation -- we obtain%
\[
E((X_{1}^{1})^{k})=\sum_{i=1}^{r}\lambda_{i}(Z_{1}^{1}(\omega_{i}))^{k},
\]
for all $0\leq k\leq m$, only if the equation holds for $k=1$, hence%
\[
E(\log(X_{1}^{1}))=\sum_{i=1}^{r}\lambda_{i}\log Z_{1}^{1}(\omega_{i})
\]
by definition of the logarithm.\label{kusuoka remark}
\end{remark}

\section{Calculation of the Greeks}

For the calculation of Greeks we proceed by methods from Malliavin Calculus
(see \cite{Nua:95} and \cite{Mal:97} for all details). We shall consider
stochastic differential equations of the type%
\begin{align*}
dY_{t}^{y}  &  =\sum_{i=0}^{d}V_{i}(Y_{t}^{y})\circ dB_{t}^{i},\\
Y_{0}^{x}  &  =y,
\end{align*}
on a stochastic basis $(\Omega,\mathcal{F},P)$, where we are given a
$d$-dimensional Brownian motion $(B_{t}^{i})_{0\leq t\leq T,i=1,\dots,d}$ in
its natural filtration, up to a finite time horizon $T>0$. For the vector
fields $V_{i}$ we shall assume the following regularity assertion:

\begin{itemize}
\item The vector fields%
\[
V_{i}:\mathbb{R}^{N}\rightarrow\mathbb{R}^{N}%
\]
are $C^{\infty}$-bounded.
\end{itemize}

If the distribution%
\begin{equation}
\mathcal{D}_{LA}(y):=\left\langle V_{1}(y),\dots,V_{d}(y),[V_{i}%
,V_{j}](y),\dots\text{ for }i,j,\dots=0,\dots,d\right\rangle \label{HC}%
\end{equation}
has constant rank $N$ for $y\in\mathbb{R}^{N}$, we say that
\emph{H\"{o}rmander's condition} holds. Here and in the sequel we apply the
notion of Lie brackets of two vector fields $V,W:\mathbb{R}^{N}\rightarrow
\mathbb{R}^{N}$, $[V,W](y):=dV(y)\cdot W(y)-dW(y)\cdot V(y)$ for
$y\in\mathbb{R}^{N}$.

\begin{theorem}
Fix $t>0$, $v\in\mathbb{R}^{N}$ and $y\in\mathbb{R}^{N}$. If condition
(\ref{HC}) holds, there is a weight $\pi\in\mathcal{D}^{\infty}$, such that
for all bounded, measurable functions $f:\mathbb{R}^{N}\rightarrow\mathbb{R}$
the equation%
\[
\frac{d}{d\epsilon}|_{\epsilon=0}E(f(Y_{t}^{y+\epsilon v})))=E(f(Y_{t}^{y}%
)\pi)
\]
holds true depending on $t,v,y$ and the whole stochastic process $(Y_{t}%
^{y})_{0\leq t\leq T}$.
\end{theorem}

\begin{proof}
For the proof see \cite{GobMun:02} and \cite{TeiTou:04}.
\end{proof}

We shall apply usual notions of Malliavin calculus (see \cite{Mal:97} and
\cite{Nua:95}). The first variation $J_{0\rightarrow t}(y)$ denotes the
derivative of $(Y_{t}^{y})_{0\leq t\leq T}$ with respect to $y$, hence%
\begin{align*}
dJ_{0\rightarrow t}(y)  &  =\sum_{i=0}^{d}dV_{i}(Y_{t}^{y})J_{0\rightarrow
t}(y)\circ dB_{t}^{i}\\
J_{0\rightarrow0}(y)  &  =id_{N},
\end{align*}
where we apply the Jacobian $dV_{i}$ of the vector field $V_{i}$ as defined in
\ref{jacobian}. This is an almost surely invertible process and we obtain the
representation%
\[
D_{s}^{k}Y_{t}^{y}=J_{0\rightarrow t}(y)(J_{0\rightarrow s}(y))^{-1}%
V_{k}(Y_{s}^{y})1_{[0,t]}(s)
\]
for $0\leq s\leq T$, $0\leq t\leq T$, $1\leq k\leq d$ of the Malliavin
derivative. We have the fundamental partial integration formula%
\[
\sum_{k=1}^{d}E(\int_{0}^{T}D_{s}^{k}Y_{t}^{y}a^{k}ds)=E(Y_{t}^{y}\delta(a)),
\]
where $(a_{s})_{0\leq s\leq T}:=(a_{0\leq s\leq T}^{k})_{k=1,\dots,d}%
\in\operatorname*{dom}(\delta)$ is a Skorohod integrable process. $\delta(a)$
is real valued random variable. Notice that if $a$ is predictable and
square-inegrable, then%
\[
\delta(a)=\sum_{k=1}^{d}\int_{0}^{T}a_{s}^{k}dB_{s}^{k},
\]
hence the Skorohod integral coincides with the Ito integral.

For the pull-back of the stochastic flow, we can provide a similar stochastic
Taylor expansion on the space of vector fields as in Section \ref{cubature},
or -- in the case of a Lie group -- on the Lie algebra. We apply the following
notation:%
\begin{align*}
\mathcal{L}_{\emptyset}V  &  :=V\\
\mathcal{L}_{i}V  &  :=[V,V_{i}],\\
\mathcal{L}_{(I)}V  &  :=\mathcal{L}_{i_{1}}\circ\cdots\circ\mathcal{L}%
_{i_{k}}(V),
\end{align*}
for $I=(i_{1},\dots,i_{k})\in\mathcal{A}$.

\begin{theorem}
Given $C^{\infty}$-bounded vector fields $V_{0},\dots,V_{d}$ on $\mathbb{R}%
^{N}$, then the diffusion process%
\begin{align*}
dY_{t}^{y}  &  =\sum_{i=0}^{d}V_{i}(Y_{t}^{y})\circ dB_{t}^{i},\\
Y_{0}^{y}  &  =y,
\end{align*}
for $y\in\mathbb{R}^{N}$, admits the following series expansion for the pull
back $J_{0\rightarrow t}(y)^{-1}V(Y_{t}^{y})$ on smooth vector fields
$V:\mathbb{R}^{N}\rightarrow\mathbb{R}^{N}$.%
\[
J_{0\rightarrow t}(y)^{-1}V(Y_{t}^{y})=\sum_{\substack{I\in\mathcal{A}%
\\\deg(I)\leq n}}\mathcal{L}_{(I)}V(y)\int_{0\leq t_{1}\leq\dots\leq t_{k}\leq
t}\circ dB_{t_{1}}^{i_{1}}\circ\cdots\circ dB_{t_{k}}^{i_{k}}+R_{n}(t,V),
\]
for $0\leq t\leq T$. For the remainder term we obtain%
\begin{gather*}
R_{n}(t,V)\\
=\sum_{\substack{(i_{1},\dots,i_{k})\in\mathcal{A}\\\deg(i_{1},\dots
,i_{k})\leq n\\(i_{0},\dots,i_{k})\in\mathcal{A}\\\deg(i_{0},\dots,i_{k}%
)>n}}\int_{0\leq t_{0}\leq\dots\leq t_{k}\leq t}J_{0\rightarrow t_{0}}%
(y)^{-1}\mathcal{L}_{(i_{0},\dots,i_{k})}V(Y_{t_{0}}^{y})\circ dB_{t_{0}%
}^{i_{0}}\circ\cdots\circ dB_{t_{k}}^{i_{k}},
\end{gather*}
for $0\leq t\leq T$.
\end{theorem}

This Taylor expansion leads on the Lie group $G_{d,1}^{m}$ to an action on the
Lie algebra $\mathfrak{g}_{d,1}^{m}$, namely%
\[
X_{t}^{1}w(X_{t}^{1})^{-1}=\sum_{\substack{I\in\mathcal{A}\\\deg(I)\leq
m-1}}[e_{i_{1}},\cdots\lbrack e_{i_{k}},w]\cdots]\int_{0\leq t_{1}\leq
\dots\leq t_{k}\leq t}\circ dB_{t_{1}}^{i_{1}}\circ\cdots\circ dB_{t_{k}%
}^{i_{k}}%
\]
by inserting the vector fields $D_{i}$ for $V_{i}$ and evaluating at $x=1$.
Notice the decrease in order from $m$ to $m-1$, since in $\mathfrak{g}%
_{d,1}^{m}$ Lie brackets only up to order $m-1$ can be non-zero. This
stochastic process on the Lie algebra $\mathfrak{g}_{d,1}^{m}$ will be applied
in the sequel in order to construct approximate, universal weights for the
calculation of the Greeks. First we need a Lemma on scaling invariance.

\begin{lemma}
The following identities in law for the stochastic process $(X_{s}^{x})_{0\leq
s\leq T}$ hold,
\begin{align*}
(\Delta_{\sqrt{t}}X_{s}^{x})_{0\leq s\leq1}  &  =_{\text{law}}(X_{st}%
^{x})_{0\leq s\leq1},\\
(\Delta_{\sqrt{t}}[(X_{s}^{e})w(X_{s}^{e})^{-1}])_{0\leq s\leq1}  &
=_{\text{law}}((\sqrt{t})^{\deg(w)}(X_{st}^{e})e_{i}(X_{st}^{e})^{-1})_{0\leq
s\leq1}%
\end{align*}
for $t>0$ and $w\in W_{n}$ (so $\deg(w)=n$).\label{scaling invariance}
\end{lemma}

\begin{proof}
By scaling invariance we know that $(\sqrt{t}B_{s}^{i})_{0\leq s\leq
1}=_{\text{law}}(B_{st}^{i})_{0\leq s\leq1}$ for $t>0$ and $i=1,\dots,d$. This
extends to all iterated stochastic integrals. Notice the importance of the
degree function, which associates degree $2$ to $e_{0}$, and hence the correct
scaling property. Notice furthermore that for $w\in\mathfrak{g}_{d,1}^{m}\cap
W_{n}$ the element $(X_{s}^{e})w(X_{s}^{e})^{-1}$ is an element of the Lie
algebra $\mathfrak{g}_{d,1}^{m}$, but not necessarily of $W_{n}$!
\end{proof}

\begin{remark}
Notice that the stochastic process $(X_{t}^{x})_{0\leq t\leq T}$ is not
hypo-elliptic, since the direction $e_{0}$ always points into directions where
the density does not admit a logarithmic derivative. Therefore we only
calculate with directions $w\in\mathfrak{g}_{d,1}^{m}/\left\langle
e_{0}\right\rangle $, which means directions with vanishing $e_{0}$ component.
In all those directions we are able to conclude a result on differentiability.
\end{remark}

\begin{remark}
In the sequel we shall construct Skorohod integrable processes $(a_{s}%
^{i})_{0\leq s\leq T}$ by defining them only on $[0,t]$ for $0<t\leq T$. It is
understood that these processes are defined on $[0,T]$, but vanish for $s>t$,
hence also $a_{s}^{i}=a_{s}^{i}1_{[0,T]}(s)$ for $0\leq s\leq T$ and
$i=1,\dots,d$.
\end{remark}

\begin{proposition}
\label{universal weight}Fix $w\in\mathfrak{g}_{d,1}^{m}/\left\langle
e_{0}\right\rangle $ (no $e_{0}$ component). For $t>0$ there are
Skorohod-integrable processes $(a_{s}^{i})_{0\leq s\leq T}$ (the
$t$-dependence is not visible in our notation) with $a_{s}^{i}=a_{s}%
^{i}1_{[0,T]}(s)$ for $0\leq s\leq T$, such that%
\[
\sum_{i=1}^{d}\int_{0}^{t}(X_{s}^{1})e_{i}(X_{s}^{1})^{-1}a_{s}^{i}%
ds=\Delta_{\sqrt{t}}(w)
\]
and%
\[
\sum_{i=1}^{d}E(\int_{0}^{t}(a_{s}^{i})^{2}ds)=\mathcal{O}(1).
\]
Furthermore we obtain that for any bounded measurable function $f:G_{d,1}%
^{m}\rightarrow\mathbb{R}$, the equation%
\[
\frac{d}{d\epsilon}|_{\epsilon=0}E(f(X_{t}^{x+\epsilon\Delta_{\sqrt{t}}%
(w)}))=E(f(X_{t}^{x})\pi_{d,1}^{m})
\]
holds true, where $\pi_{d,1}^{m}:=\delta(s\mapsto a_{s})$ denotes the Skorohod
integral of the strategies $a$.
\end{proposition}

\begin{proof}
The stochastic process $(X_{t}^{x})_{0\leq t\leq T}$ is not hypo-elliptic, but
we prove that derivatives in all directions $w\in\mathfrak{g}_{d,1}^{m}$ with
vanishing $e_{0}$ component exist. We can define the processes $a_{s}^{i}$ by
the standard method of Malliavin Calculus: we choose a scalar product
$\left\langle ,\right\rangle $ on $\mathfrak{g}_{d,1}^{m}$, which respects the
grading, i.e. $W_{k}$ is orthogonal to $W_{l}$ for $k\neq l$, and define the
reduced covariance operator $C^{t}$ as symmetric operator $C^{t}%
:\mathfrak{g}_{d,1}^{m}\rightarrow\mathfrak{g}_{d,1}^{m}$by%
\[
\left\langle y,C^{t}y\right\rangle :=\sum_{i=1}^{d}\int_{0}^{t}\left\langle
(X_{s}^{1})e_{i}(X_{s}^{1})^{-1},y\right\rangle ^{2}ds.
\]
This operator is almost surely invertible on $\mathfrak{g}_{d,1}%
^{m}/\left\langle e_{0}\right\rangle $ with $p$-integrable inverse (see
\cite{Nua:95}). Obviously the image of $C^{t}$ lies in $\mathfrak{g}_{d,1}%
^{m}/\left\langle e_{0}\right\rangle $, since there is no $e_{0}$-component in
$(X_{s}^{1})^{-1}e_{i}X_{s}^{1}$, the kernel of $C^{t}$ can be calculated by
the classical method from \cite{Nua:95} and vanishes on $\mathfrak{g}%
_{d,1}^{m}/\left\langle e_{0}\right\rangle $, hence $C^{t}|_{\mathfrak{g}%
_{d,1}^{m}/\left\langle e_{0}\right\rangle }$ is invertible with
$p$-integrable inverse, since the Norris-Lemma applies.

Furthermore we observe the following scaling property: notice that
$\Delta_{\sqrt{t}}$ is a self-adjoint operator on $\mathfrak{g}_{d,1}^{m}$ for
$t>0$. We can compare $\Delta_{\sqrt{t}}C^{1}\Delta_{\sqrt{t}}$ with $C^{t}$
in law via Lemma \ref{scaling invariance}%
\begin{align*}
\left\langle y,C^{t}y\right\rangle  &  =\sum_{i=1}^{d}\int_{0}^{t}\left\langle
(X_{s}^{1})e_{i}(X_{s}^{1})^{-1},y\right\rangle ^{2}ds\\
\left\langle \Delta_{\sqrt{t}}y,C^{1}\Delta_{\sqrt{t}}y\right\rangle  &
=\sum_{i=1}^{d}\int_{0}^{1}\left\langle (X_{s}^{1})e_{i}(X_{s}^{1}%
)^{-1},\Delta_{\sqrt{t}}y\right\rangle ^{2}ds\\
&  =\sum_{i=1}^{d}\int_{0}^{1}\left\langle \Delta_{\sqrt{t}}[(X_{s}^{1}%
)e_{i}(X_{s}^{1})^{-1}],y\right\rangle ^{2}ds\\
&  =_{\text{law}}\sum_{i=1}^{d}\int_{0}^{1}\left\langle \sqrt{t}(X_{st}%
^{1})e_{i}(X_{st}^{1})^{-1},y\right\rangle ^{2}ds\\
&  =\sum_{i=1}^{d}\int_{0}^{t}\left\langle (X_{u}^{1})e_{i}(X_{u}^{1}%
)^{-1},y\right\rangle ^{2}du.
\end{align*}
So we have that
\[
((X_{s}^{1})e_{i}(X_{s}^{1}))_{0\leq s\leq t},C^{t})=_{\text{law}}((\frac
{1}{\sqrt{t}}\Delta_{\sqrt{t}}[(X_{\frac{s}{t}}^{1})e_{i}(X_{\frac{s}{t}}%
^{1})^{-1}])_{0\leq s\leq t},\Delta_{\sqrt{t}}C^{1}\Delta_{\sqrt{t}})
\]
for $t>0$. We define the (non-adapted) strategies via%
\[
a_{s}^{i}:=\left\langle (X_{s}^{1})e_{i}(X_{s}^{1}),(C^{t}|_{\mathfrak{g}%
_{d,1}^{m}/\left\langle e_{0}\right\rangle })^{-1}\Delta_{\sqrt{t}%
}(w)\right\rangle
\]
for $0\leq s\leq t$, $i=1,\dots,d$, and obtain%
\begin{gather*}
\sum_{i=1}^{d}\int_{0}^{t}\left\langle (X_{s}^{1})e_{i}(X_{s}^{1}%
)^{-1},v\right\rangle a_{s}^{i}ds=\\
\sum_{i=1}^{d}\int_{0}^{t}\left\langle (X_{s}^{1})e_{i}(X_{s}^{1}%
)^{-1},v\right\rangle \left\langle (X_{s}^{1})e_{i}(X_{s}^{1})^{-1}%
,(C^{t})^{-1}\Delta_{\sqrt{t}}(w)\right\rangle ds=\\
\left\langle v,C^{t}(C^{t})^{-1}\Delta_{\sqrt{t}}(w)\right\rangle
=\left\langle v,\Delta_{\sqrt{t}}(w)\right\rangle
\end{gather*}
for $v\in\mathfrak{g}_{d,1}^{m}$. Hence the strategies satisfy the desired
equation by duality on $\mathfrak{g}_{d,1}^{m}$. Concerning the order estimate
for the strategies we observe that%
\begin{align*}
&  \left\langle (X_{s}^{1})e_{i}(X_{s}^{1})^{-1},(C^{t})^{-1}\Delta_{\sqrt{t}%
}(w)\right\rangle \\
&  =_{\text{law}}\left\langle \frac{1}{\sqrt{t}}\Delta_{\sqrt{t}}[(X_{\frac
{s}{t}}^{1})e_{i}(X_{\frac{s}{t}}^{1})^{-1}],\Delta_{\sqrt{t}}^{-1}%
(C^{1})^{-1}\Delta_{\sqrt{t}}^{-1}\Delta_{\sqrt{t}}(w)\right\rangle \\
&  =\left\langle \frac{1}{\sqrt{t}}(X_{\frac{s}{t}}^{1})e_{i}(X_{\frac{s}{t}%
}^{1})^{-1},(C^{1})^{-1}w\right\rangle ,
\end{align*}
so consequently%
\begin{align*}
\sum_{i=1}^{d}E(\int_{0}^{t}(a_{s}^{i})^{2}ds)  &  =_{\text{law}}\sum
_{i=1}^{d}E(\int_{0}^{t}\left\langle \frac{1}{\sqrt{t}}(X_{\frac{s}{t}}%
^{1})e_{i}(X_{\frac{s}{t}}^{1})^{-1},(C^{1})^{-1}w\right\rangle ^{2}ds\\
&  =\sum_{i=1}^{d}E(\int_{0}^{t}\left\langle \frac{1}{\sqrt{t}}(X_{\frac{s}%
{t}}^{1})e_{i}(X_{\frac{s}{t}}^{1})^{-1},(C^{1})^{-1}w\right\rangle
^{2}td(\frac{s}{t})\\
&  =\sum_{i=1}^{d}E(\int_{0}^{1}\left\langle (X_{u}^{1})e_{i}(X_{u}^{1}%
)^{-1},(C^{1})^{-1}w\right\rangle ^{2}du
\end{align*}

\end{proof}

\begin{remark}
The strategies are \textbf{universal} in the sense that they only depend on
the Brownian motion, the choice of a scalar product on $\mathfrak{g}_{d,1}%
^{m}$, and -- in a linear way -- on the \textquotedblright
input\textquotedblright\ direction $w$. We denote the Skorohod integrals of
those strategies by $\pi_{d,1}^{m}:=\delta(a)$ and call them \emph{universal
weights}. Notice also that $E((\pi_{d,1}^{m})^{2})=\mathcal{O}(1)$, which is
seen by a similar calculation as the previous one. Hence there is a universal
time dependence, which is achieved by calculating the strategies with respect
to the direction $\Delta_{\sqrt{t}}w$ instead of $w$. We omit the dependence
on $w$ in our notation.\label{universal remark}
\end{remark}

\begin{example}
We consider the simplest non-trivial example: $m=2$ and $d=2$. Hence the
nilpotent algebra is generated by%
\begin{align*}
W_{0}  &  =\left\langle 1\right\rangle ,\\
W_{1}  &  =\left\langle e_{1},e_{2}\right\rangle ,\\
W_{2}  &  =\left\langle e_{1}e_{2},e_{2}e_{1},e_{0}\right\rangle ,
\end{align*}
and the Lie algebra by $e_{1},e_{2},e_{0},[e_{1},e_{2}]$. The dimensions are
$6$ and $4$ respectively. We can solve explicitly%
\[
dX_{t}^{x}=\sum_{i=0}^{2}X_{t}^{x}e_{i}\circ dB_{t}^{i}%
\]
by%
\begin{align*}
X_{t}^{x}  &  =x+xe_{1}B_{t}^{1}+xe_{2}B_{t}^{2}+xe_{0}t+xe_{1}e_{2}\int
_{0}^{t}\int_{0}^{t_{2}}\circ dB_{t_{1}}^{1}\circ dB_{t_{2}}^{2}+\\
&  +xe_{2}e_{1}\int_{0}^{t}\int_{0}^{t_{2}}\circ dB_{t_{1}}^{2}\circ
dB_{t_{2}}^{1}+\frac{1}{2}xe_{1}^{2}(B_{t}^{1})^{2}+\frac{1}{2}xe_{2}%
^{2}(B_{t}^{2})^{2}.
\end{align*}
The pull-back yields a shorter expression, namely%
\begin{align*}
(X_{t}^{1})e_{1}(X_{t}^{1})^{-1}  &  =e_{1}+[e_{2},e_{1}]B_{t}^{2},\\
(X_{t}^{1})e_{2}(X_{t}^{1})^{-1}  &  =e_{2}+[e_{1},e_{2}]B_{t}^{1},
\end{align*}
for $0\leq t\leq T$, which can also be calculated directly. Consequently the
above equation reduces to find a Skorohod integrable strategy such that%
\begin{align*}
\int_{0}^{t}a_{s}^{1}ds  &  =\sqrt{t}w_{1}\\
\int_{0}^{t}a_{s}^{2}ds  &  =\sqrt{t}w_{2}\\
\int_{0}^{t}(a_{s}^{1}B_{s}^{2}-a_{s}^{2}B_{s}^{1})ds  &  =tw_{3}%
\end{align*}
for $w=w_{1}e_{1}+w_{2}e_{2}+w_{3}[e_{1},e_{2}].$ This can be done by
introducing a scalar product on $\mathfrak{g}_{d,1}^{m}$ such that
$e_{1},e_{2},[e_{1},e_{2}],e_{0}$ is an orthonormal basis. We then obtain the
symmetric matrix%
\[
C^{t}=\left(
\begin{array}
[c]{cccc}%
t & 0 & \int_{0}^{t}B_{s}^{2}ds & 0\\
0 & t & -\int_{0}^{t}B_{s}^{1}ds & 0\\
\int_{0}^{t}B_{s}^{2}ds & -\int_{0}^{t}B_{s}^{1}ds & \int_{0}^{t}((B_{s}%
^{1})^{2}+(B_{s}^{2})^{2})ds & 0\\
0 & 0 & 0 & 0
\end{array}
\right)
\]
with respect to the given basis. So $C^{t} $ is almost surely invertible on
$\left\langle e_{1},e_{2},[e_{1},e_{2}]\right\rangle $ -- as claimed in
Proposition \ref{universal weight}, indeed the determinant is given through%
\[
\det(C^{t}|_{\left\langle e_{1},e_{2},[e_{1},e_{2}]\right\rangle })=t^{2}%
\int_{0}^{t}((B_{s}^{1})^{2}+(B_{s}^{2})^{2})ds-t(\int_{0}^{t}B_{s}^{1}%
ds)^{2}-t(\int_{0}^{t}B_{s}^{2}ds)^{2},
\]
which is positive if and only if one of the two Brownian motions is not
vanishing identically on $[0,t]$. The order assertion in $t$ is also nicely
visible.\label{example m=d=2}
\end{example}

We consider the first derivative of the function $y\mapsto E(f(Y_{t}^{y}))$
for $t>0$. By H\"{o}rmander's theorem we know that this function is smooth for
all bounded measurable functions $f:\mathbb{R}^{N}\rightarrow\mathbb{R}$ (see
for instance \cite{Mal:97} or \cite{Nua:95}). We generalize this assertion on
the one hand, since we leave away H\"{o}rmander's condition. On the other hand
we only prove an approximative result with a certain time asymptotics.

\begin{theorem}
Fix $t>0$, $v\in\mathbb{R}^{N}$ and $y\in\mathbb{R}^{N}$. Furthermore we fix
an order $m\geq1$ of approximation. We assume that the vector $v$ is a linear
combination of Lie brackets up to order $m-1$ (except the direction $V_{0}$),
i.e.%
\[
v=\sum_{\substack{I\in\mathcal{A}\setminus\{\emptyset,(0)\}\\\deg(I)\leq
m-1}}t^{\frac{\deg(I)}{2}}w_{I}[V_{i_{1}},[V_{i_{2}},[\cdots,V_{i_{k}}%
]\cdots](y),
\]
and define%
\[
\Delta_{\sqrt{t}}(w):=\sum_{\substack{I\in\mathcal{A}\setminus\{\emptyset
,(0)\}\\\deg(I)\leq m-1}}t^{\frac{\deg(I)}{2}}w_{I}[e_{i_{1}},[e_{i_{2}%
},[\cdots,e_{i_{k}}]\cdots].
\]
Then there is universal weight $\pi=\pi_{d,1}^{m}\in\mathcal{D}^{\infty}$
associated to $\Delta_{\sqrt{t}}(w)$, such that for all $C^{\infty}$-bounded
functions $f:\mathbb{R}^{N}\rightarrow\mathbb{R}$ the equation%
\[
\frac{d}{d\epsilon}|_{\epsilon=0}E(f(Y_{t}^{y+\epsilon v}))=E(f(Y_{t}^{y}%
)\pi_{d,1}^{m})+\mathcal{O}(t^{\frac{m+1}{2}})
\]
holds true, where the constant in the order estimate only depends on the first
derivative of $f$.\label{existence-weights}
\end{theorem}

\begin{proof}
Assume $f$ is $C_{b}^{1}$, then%
\[
\frac{d}{d\epsilon}|_{\epsilon=0}E(f(Y_{t}^{y+\epsilon v}))=E(df(Y_{t}%
^{y})\cdot J_{0\rightarrow t}(y)\cdot v),
\]
where $df$ denotes the differential of the function $f$. Instead of solving
the equation%
\[
v=\sum_{k=1}^{d}\int_{0}^{t}J_{0\rightarrow s}(y)^{-1}V_{k}(Y_{s}^{y}%
)a_{s}^{k}ds,
\]
we take the universal weight $\pi_{d,1}^{m}=\delta(s\mapsto a_{s})$ of
Proposition \ref{universal weight} and solve the equation%
\[
v=\sum_{k=1}^{d}\int_{0}^{t}(J_{0\rightarrow s}(y)^{-1}V_{k}(Y_{s}^{y}%
)-R_{m}^{k}(s,V_{k}))a_{s}^{k}ds
\]
by the above universal construction. More precisely, we choose real numbers
$w_{I}$ for $I\in\mathcal{A}\setminus\{\emptyset\}$ such that%
\[
v=\sum_{\substack{I\in\mathcal{A}\setminus\{\emptyset,(0)\}\\\deg(I)\leq
m-1}}t^{\frac{\deg(I)}{2}}w_{I}[V_{i_{1}},[V_{i_{2}},[\cdots,V_{i_{k}}%
]\cdots](y)
\]
and define%
\[
\Delta_{\sqrt{t}}w:=\sum_{\substack{I\in\mathcal{A}\setminus\{\emptyset
,(0)\}\\\deg(I)\leq m-1}}t^{\frac{\deg(I)}{2}}w_{I}[e_{i_{1}},[e_{i_{2}%
},[\cdots,e_{i_{k}}]\cdots].
\]
Then we take the universal weight $\pi:=\pi_{d,1}^{m}$ associated to $w$ and
solve consequently the above equation. This leads finally to%
\begin{gather*}
\frac{d}{d\epsilon}|_{\epsilon=0}E(f(Y_{t}^{y+\epsilon v}))=E(df(Y_{t}%
^{y})\cdot J_{0\rightarrow t}(y)\cdot w)\\
=E(df(Y_{t}^{y})\cdot J_{0\rightarrow t}(y)\cdot\sum_{k=1}^{d}\int_{0}%
^{t}J_{0\rightarrow s}(y)^{-1}V_{k}(Y_{s}^{y})a_{s}^{k}ds)+\mathcal{O}%
(t^{\frac{m+1}{2}})\\
=E(\sum_{k=1}^{d}\int_{0}^{t}df(Y_{t}^{y})\cdot J_{0\rightarrow t}(y)\cdot
J_{0\rightarrow s}(y)^{-1}V_{k}(Y_{s}^{y})a_{s}^{k}ds)+\mathcal{O}%
(t^{\frac{m+1}{2}})\\
=E(\sum_{k=1}^{d}\int_{0}^{t}D_{s}^{k}(Y_{t}^{y})a_{s}^{k}ds)+\mathcal{O}%
(t^{\frac{m+1}{2}})\\
=E(f(Y_{t}^{y})\pi)+\mathcal{O}(t^{\frac{m+1}{2}}),
\end{gather*}
hence the result, since due to Proposition \ref{universal weight} we obtain --
by the Cauchy-Schwartz inequality and integration with respect to $t$ -- an
order estimate of the type $\mathcal{O}(t^{\frac{m+1}{2}})$: indeed, the
remainder satisfies%
\begin{align*}
&  \left\vert E(\sum_{k=1}^{d}\int_{0}^{t}df(Y_{t}^{y})\cdot J_{0\rightarrow
t}(y)\cdot R_{m}^{k}(s,V_{k})a_{s}^{k}ds)\right\vert \\
&  \leq(\sum_{k=1}^{d}E(\int_{0}^{t}(df(Y_{t}^{y})\cdot J_{0\rightarrow
t}(y)\cdot R_{m}^{k}(s,V_{k}))^{2}ds)^{\frac{1}{2}}(\sum_{i=1}^{d}E(\int
_{0}^{t}(a_{s}^{i})^{2}ds))^{\frac{1}{2}}\\
&  \leq\sum_{k=1}^{d}M(\int_{0}^{t}s^{m}ds)^{\frac{1}{2}}=\mathcal{O}%
(t^{\frac{m+1}{2}}).
\end{align*}

\end{proof}

\section{Cubature formulas for Greeks}

First we have to get expressions for derivatives of heat equation evolutions
on free, nilpotent Lie groups, which is an easy task for linear functions by
the considerations of Section \ref{cubature}.

\begin{proposition}
Fix a linear function $\lambda$ on $\mathbb{A}_{d,1}^{m}$, then%
\[
\frac{d}{d\epsilon}|_{\epsilon=0}E(\lambda(X_{t}^{x+\epsilon w_{x}}%
))=\lambda(w_{x}\exp(t(e_{0}+\frac{1}{2}\sum_{i=1}^{d}e_{i}^{2})))
\]
for $w_{x}\in x(\mathfrak{g}_{d,1}^{m}/\left\langle e_{0}\right\rangle )$.
Consequently%
\begin{equation}
E(X_{t}^{1}\pi_{d,1}^{m})=\Delta_{\sqrt{t}}(w)\exp(t(e_{0}+\frac{1}{2}%
\sum_{i=1}^{d}e_{i}^{2})) \label{expectationformula der}%
\end{equation}
for $w\in\mathfrak{g}_{d,1}^{m}/\left\langle e_{0}\right\rangle $%
.\label{derivative}
\end{proposition}

\begin{proof}
We know that the function $x\mapsto E(\lambda(X_{t}^{x}))$ solves the
associated heat equation. The solution of the heat equation on $G_{d,1}^{m}$
with initial value $\lambda$ is given by%
\[
x\mapsto\lambda(x\exp(t(e_{0}+\frac{1}{2}\sum_{i=1}^{d}e_{i}^{2}))),
\]
hence the derivative in direction $w_{x}$ can be calculated and yields%
\[
\lambda(w_{x}\exp(t(e_{0}+\frac{1}{2}\sum_{i=1}^{d}e_{i}^{2}))),
\]
wherefrom the result follows by duality and Theorem \ref{existence-weights}.
\end{proof}

Next we apply Tchakaloff's theorem twice in order to obtain the appropriate
cubature result.

\begin{theorem}
Fix a free, nilpotent Lie group $G_{d,1}^{m}$ and $w\in\mathfrak{g}_{d,1}%
^{m}/\left\langle e_{0}\right\rangle $, then there are points $x_{1}%
,\dots,x_{r}\in G_{d,1}^{m}$ and weights $\mu_{1},\dots,\mu_{r}\neq0$ such
that%
\[
\Delta_{\sqrt{t}}(w)\exp(t(e_{0}+\frac{1}{2}\sum_{i=1}^{d}e_{i}^{2}%
))=\sum_{j=1}^{r}x_{j}\mu_{j}.
\]
Furthermore $r\leq2\dim_{\mathbb{R}}\mathbb{A}_{d,1}^{m}$.
\end{theorem}

\begin{proof}
We write $\pi=\pi_{d,1}^{m}$. We define two probability measures on Wiener
space, absolutely continuous with respect to Wiener measure $P$,%
\begin{align*}
\frac{dQ_{+}}{dP}  &  =\frac{1}{E(\pi_{+})}\pi_{+},\\
\frac{dQ_{-}}{dP}  &  =\frac{1}{E(\pi_{-})}\pi_{-},
\end{align*}
if the respective positive and negative parts have non-vanishing expectation.
Then%
\[
E_{P}(X_{t}^{e}\pi)=E(\pi_{+})E_{Q_{+}}(X_{t}^{e})-E(\pi_{-})E_{Q_{-}}%
(X_{t}^{e}).
\]
There is a null set $N$ on Wiener space, such that on $N^{c}$ the process
$X_{t}^{e}\in G_{d,1}^{m}$, hence also $X_{t}^{e}\in G_{d,1}^{m}$ almost
surely with respect to $Q_{\pm}$. Consequently by Tchakaloff's theorem we find
points $x_{1}^{\pm},\dots,x_{r}^{\pm}$ and weights $\lambda_{j}^{\pm}$ such
that%
\[
E_{P}(X_{t}^{e}\pi)=\sum_{j=1}^{r}(E(\pi_{+}^{(1)})x_{j}^{+}\lambda_{j}%
^{+}-E(\pi_{-}^{(1)})x_{j}^{-}\lambda_{j}^{-}),
\]
which yields the desired formula by Proposition \ref{derivative}.
\end{proof}

\begin{remark}
We can additionally require that $\sum_{i=1}^{r}\mu_{i}=0$, $\sum_{i=1}%
^{r}|\mu_{i}|\leq2$ and $|\mu_{i}|\leq1$ for $i=1,\dots,r$.
\end{remark}

Hence we can formulate a Cubature result for the calculation of Greeks.

\begin{definition}
Fix $m\geq1$ and $0<t\leq T$, a \emph{cubature formula for the first
derivative} in direction $\Delta_{\sqrt{t}}w\in\mathfrak{g}_{d,1}%
^{m}/\left\langle e_{0}\right\rangle $ is given by a finite number of points
$x_{1},\dots,x_{r}\in G_{d,1}^{m}$ and finitely many weights $\mu_{1}%
,\dots,\mu_{r}\neq0$, such that%
\[
E(X_{t}^{1}\pi_{d,1}^{m})=\sum_{j=1}^{r}\mu_{j}x_{j},
\]
or equivalently due to formula (\ref{expectationformula der})%
\[
\Delta_{\sqrt{t}}(w)\exp(t(e_{0}+\frac{1}{2}\sum_{i=1}^{d}e_{i}^{2}%
))=\sum_{j=1}^{r}\mu_{j}x_{j}.
\]
Again we use trajectories $\omega_{i}$ to represent the points $x_{i}$ as
endpoints of evolutions of ordinary differential equations $Z_{t}(\omega
_{i})=x_{i}$ for $i=1,\dots,r$. Furthermore $r\leq2\dim_{\mathbb{R}}%
\mathbb{A}_{d,1}^{m}$.
\end{definition}

\begin{theorem}
\label{cubature for greeks}Given $C^{\infty}$-bounded vector fields
$V_{0},\dots,V_{d}$ on $\mathbb{R}^{N}$, then the diffusion process%
\[
dY_{t}^{y}=\sum_{i=0}^{d}V^{i}(Y_{t}^{y})\circ dB_{t}^{i}%
\]
for $y\in\mathbb{R}^{N}$ admits the following cubature formulas for
derivatives in direction $v\in\mathbb{R}^{N}$. Fix $m\geq1$ and assume that%
\begin{equation}
v=\sum_{\substack{I\in\mathcal{A}\setminus(\emptyset,(0))\\\deg(I)\leq
m-1}}t^{\frac{\deg(I)}{2}}w_{I}[V_{i_{1}},[V_{i_{2}},[\cdots,V_{i_{k}}%
]\cdots](y), \label{decompo}%
\end{equation}
then we obtain
\[
\frac{d}{d\epsilon}|_{\epsilon=0}E(f(Y_{t}^{y+\epsilon v}))=\sum_{j=1}^{r}%
\mu_{j}f(Y_{t}^{y}(\omega_{j}))+\mathcal{O}(t^{\frac{m+1}{2}}),
\]
taking a cubature formula, which was derived in $G_{d,1}^{m}$ and $C^{\infty}%
$-bounded $f$. The constant depends in general on derivatives of $f$ up to
order $m+1$. We can determine the weights $\mu_{j}$ and the trajectories
$\omega_{j}$ by%
\begin{gather*}
\Delta_{\sqrt{t}}(w)\exp(t(e_{0}+\frac{1}{2}\sum_{i=1}^{d}e_{i}^{2}))\\
=\sum_{i=1}^{r}\mu_{j}\sum_{\substack{I\in\mathcal{A}\setminus(\emptyset
,(0))\\\deg(I)\leq m}}e_{i_{1}}\cdots e_{i_{k}}\int_{0\leq t_{1}\leq\dots\leq
t_{k}\leq t}d\omega_{j}^{i_{1}}(t_{1})\cdots d\omega_{j}^{i_{k}}(t_{k}),\\
\Delta_{\sqrt{t}}(w)=\sum_{\substack{I\in\mathcal{A}\setminus(\emptyset
,(0))\\\deg(I)\leq m-1}}t^{\frac{\deg(I)}{2}}w_{I}[e_{i_{1}},[e_{i_{2}%
},[\cdots,e_{i_{k}}]\cdots].
\end{gather*}

\end{theorem}

\begin{remark}
Since we shall in practise apply this procedure to smooth functions -- as
discussed in the recipe after formulas (\ref{recipe1}) and (\ref{recipe2}) --
this already yields the interesting result for applications.
\end{remark}

\begin{proof}
By the previous constructions we obtain%
\[
E(X_{t}^{x}\pi_{d,1}^{m})=\sum_{i=1}^{r}\mu_{j}\sum_{\substack{I\in
\mathcal{A}\\\deg(I)\leq m}}e_{i_{1}}\cdots e_{i_{k}}\int_{0\leq t_{1}%
\leq\dots\leq t_{k}\leq t}d\omega_{j}^{i_{1}}(t_{1})\cdots d\omega_{j}^{i_{k}%
}(t_{k}).
\]
We then insert instead of $e_{i}$ the vector fields $V_{i}$ and apply those
vector fields to the function $f$ at $y\in\mathbb{R}^{N}$. Then we obtain on
the left hand side -- due to stochastic Taylor expansion -- and on the right
hand side due to Taylor expansion of the non-autonomous equation,%
\[
E((f(Y_{t}^{y})-R_{m}(t,f))\pi_{d,1}^{m})=\sum_{i=1}^{r}\mu_{j}(f(Y_{t}%
^{y}(\omega_{j}))-R_{m,i}(t,f)).
\]
The left hand side yields then the order estimate by partial integration and
Proposition \ref{universal weight}, the right hand side by the fact that
iterated integrals with respect to $\omega_{j}$ behave like $t^{\frac{m+1}{2}%
}$.
\end{proof}

We additionally have an assertion on the construction of cubature paths due to
the scaling properties. Assume that we found trajectories $\omega_{i}$ and
weights $\mu_{i}$ for $m\geq1$, $d\geq1$ and $t=1$ fixed, then we can
construct solutions for all $t>0$ with the same $m,d$: the equation%
\[
w\exp(e_{0}+\frac{1}{2}\sum_{i=1}^{d}e_{i}^{2})=\sum_{i=1}^{r}\mu_{j}%
\sum_{\substack{I\in\mathcal{A}\\\deg(I)\leq m}}e_{i_{1}}\cdots e_{i_{k}}%
\int_{0\leq t_{1}\leq\dots\leq t_{k}\leq1}d\omega_{j}^{i_{1}}(t_{1})\cdots
d\omega_{j}^{i_{k}}(t_{k})
\]
holds, hence the trajectories%
\begin{align*}
\eta_{j}^{0}(st)  &  :=t\omega_{j}^{0}(s),\\
\eta_{j}^{i}(st)  &  :=\sqrt{t}\omega_{j}^{i}(s),
\end{align*}
for $i=1,\dots,d$, $j=1,\dots,r$ and $0\leq s\leq1$ satisfy%
\begin{gather*}
\Delta_{\sqrt{t}}(w)\exp(t(e_{0}+\frac{1}{2}\sum_{i=1}^{d}e_{i}^{2}))\\
=\sum_{i=1}^{r}\mu_{j}\sum_{\substack{I\in\mathcal{A}\\\deg(I)\leq m}%
}e_{i_{1}}\cdots e_{i_{k}}\int_{0\leq t_{1}\leq\dots\leq t_{k}\leq t}d\eta
_{j}^{i_{1}}(t_{1})\cdots d\eta_{j}^{i_{k}}(t_{k}),
\end{gather*}
since%
\begin{gather*}
\sum_{i=1}^{r}\mu_{j}\sum_{\substack{I\in\mathcal{A}\\\deg(I)\leq m}}e_{i_{1}%
}\cdots e_{i_{k}}\int_{0\leq t_{1}\leq\dots\leq t_{k}\leq t}d\eta_{j}^{i_{1}%
}(t_{1})\cdots d\eta_{j}^{i_{k}}(t_{k})\\
=\sum_{i=1}^{r}\mu_{j}\sum_{\substack{I\in\mathcal{A}\\\deg(I)\leq m}%
}\Delta_{\sqrt{t}}(e_{i_{1}}\cdots e_{i_{k}})\int_{0\leq t_{1}\leq\dots\leq
t_{k}\leq1}d\omega_{j}^{i_{1}}(t_{1})\cdots d\omega_{j}^{i_{k}}(t_{k})\\
=\Delta_{\sqrt{t}}(\sum_{i=1}^{r}\mu_{j}\sum_{\substack{I\in\mathcal{A}%
\\\deg(I)\leq m}}e_{i_{1}}\cdots e_{i_{k}}\int_{0\leq t_{1}\leq\dots\leq
t_{k}\leq1}d\omega_{j}^{i_{1}}(t_{1})\cdots d\omega_{j}^{i_{k}}(t_{k}))\\
=\Delta_{\sqrt{t}}(w\exp(e_{0}+\frac{1}{2}\sum_{i=1}^{d}e_{i}^{2}))\\
=\Delta_{\sqrt{t}}(w)\exp(t(e_{0}+\frac{1}{2}\sum_{i=1}^{d}e_{i}^{2})).
\end{gather*}
This leads to the desired assertion.

\section{Applications}

\subsection{Example for $m=1$ and $d\geq1$}

In \cite{FouLasLebLioTou:99} the authors provide the following expression for
the Malliavin weight $\pi$: given $C^{\infty}$-bounded vector fields
$V_{1},\dots,V_{d}$ on $\mathbb{R}^{d}$ such that a uniform ellipticity
condition holds, i.e. there is $\delta>0$ such that%
\[
\sum_{i=1}^{d}\left\langle \xi,V_{i}\right\rangle ^{2}\geq\delta\left\langle
\xi,\xi\right\rangle
\]
for all $\xi\in\mathbb{R}^{d}$, with respect to a standard scalar product on
$\mathbb{R}^{d}$, the formula%
\begin{equation}
\frac{d}{d\epsilon}|_{\epsilon=0}E(f(Y_{t}^{y+\epsilon v}))=E(f(Y_{t}%
^{y})\frac{1}{t}\int_{0}^{t}(\sigma^{-1}(Y_{s}^{y})(J_{0\rightarrow s}(y)\cdot
v))^{T}dB_{s}), \label{precise elliptic formula}%
\end{equation}
where we understand $dB_{s}$ as column vector, where the random matrix
$\sigma^{-1}(Y_{s}^{y})$ acts on. $\sigma$ is defined via%
\[
\sigma(y):=(V_{1}(y),\dots,V_{d}(y)).
\]
Notice that this formula provides adapted strategies, therefore the Skorohod
integral can be replaced by the Ito integral. We shall compare this formula to
the formula obtained in our setting. We choose $m=1$, we can calculate
$\pi_{d,1}^{1}$ explicitly, namely%
\[
\pi_{d,1}^{1}=\sum_{i=1}^{d}B_{t}^{i}\frac{w_{i}}{\sqrt{t}}%
\]
for $i=1,\dots,d$ and $w$ has $e_{i}$-component $w_{i}$ for $i=1,\dots,d$.
Hence we approximate
\[
\frac{d}{d\epsilon}|_{\epsilon=0}E(f(Y_{t}^{y+\epsilon\sqrt{t}v}%
))=E(f(Y_{t}^{y})\sum_{i=1}^{d}B_{t}^{i}\frac{w_{i}}{\sqrt{t}})+\mathcal{O}%
(t),
\]
where $w_{i}$ stem from the solution of the equation%
\[
v=\sum_{i=1}^{d}w_{i}V_{i}(y),
\]
consequently $w=\sigma^{-1}(y)\cdot v$. Hence%
\begin{align*}
\frac{d}{d\epsilon}|_{\epsilon=0}E(f(Y_{t}^{y+\epsilon\sqrt{t}v}))  &
=E(f(Y_{t}^{y})\sum_{i=1}^{d}B_{t}^{i}\frac{\sigma^{-1}(y)\cdot v}{\sqrt{t}%
})+\mathcal{O}(t),\\
\frac{d}{d\epsilon}|_{\epsilon=0}E(f(Y_{t}^{y+\epsilon v}))  &  =E(f(Y_{t}%
^{y})\sum_{i=1}^{d}B_{t}^{i}\frac{\sigma^{-1}(y)\cdot v}{t})+\mathcal{O}%
(\sqrt{t}),
\end{align*}
which is obviously the expansion of the precise formula
(\ref{precise elliptic formula}), since%
\begin{gather*}
E(f(Y_{t}^{y})\frac{1}{t}\int_{0}^{t}(\sigma^{-1}(Y_{s}^{y})(J_{0\rightarrow
s}(y)\cdot v))^{T}dB_{s})=E(f(Y_{t}^{y})\sum_{i=1}^{d}B_{t}^{i}\frac
{\sigma^{-1}(y)\cdot v}{t})+\\
+E(f(Y_{t}^{y})\frac{1}{t}\int_{0}^{t}(\sigma^{-1}(Y_{s}^{y})(J_{0\rightarrow
s}(y)\cdot v)-\sigma^{-1}(y)\cdot v)^{T}dB_{s})\\
=E(f(Y_{t}^{y})\sum_{i=1}^{d}B_{t}^{i}\frac{\sigma^{-1}(y)\cdot v}{t})+\\
+E(\frac{1}{t}\int_{0}^{t}D_{s}f(Y_{t}^{y})(\sigma^{-1}(Y_{s}^{y}%
)(J_{0\rightarrow s}(y)\cdot v)-\sigma^{-1}(y)\cdot v)^{T}ds)\\
=E(f(Y_{t}^{y})\sum_{i=1}^{d}B_{t}^{i}\frac{\sigma^{-1}(y)\cdot v}%
{t})+\mathcal{O}(\sqrt{t})
\end{gather*}
for $f$ a $C_{b}^{\infty}$-function.

\subsection{Example for $m=2$ and $d=2$}

In order to demonstrate the method we calculate one example of the above
methodology. First we fix $m=d=2$, hence we can include the calculations of
Example \ref{example m=d=2}. We first fix a direction $w=e_{1}$. Then we find
one solution of%
\[
\sqrt{t}e_{1}=\sqrt{t}e_{1}\exp(t(e_{0}+\frac{1}{2}\sum_{i=1}^{2}e_{i}%
^{2}))=\sum_{i=1}^{r}\mu_{i}x_{i},
\]
which can be given through $r=2$, $\mu_{1}=\frac{1}{2}$, $\mu_{2}=-\frac{1}%
{2}$ and $x_{1}=\exp(\sqrt{t}e_{1})$, $x_{2}=\exp(-\sqrt{t}e_{1})$.
Consequently the trajectories are given by%
\begin{align*}
\omega_{1}(s)  &  =\left(
\begin{array}
[c]{c}%
0\\
\frac{s}{\sqrt{t}}\\
0
\end{array}
\right)  ,\\
\omega_{2}(s)  &  =\left(
\begin{array}
[c]{c}%
0\\
-\frac{s}{\sqrt{t}}\\
0
\end{array}
\right)  ,
\end{align*}
which leads to%
\begin{align*}
Z_{t}^{1}(\omega_{1})  &  =1+\sqrt{t}e_{1}+\frac{t}{2}e_{1}^{2}=\exp(\sqrt
{t}e_{1}),\\
Z_{t}^{1}(\omega_{2})  &  =1-\sqrt{t}e_{1}+\frac{t}{2}e_{1}^{2}=\exp(-\sqrt
{t}e_{1}).
\end{align*}
This in turn leads to the following approximation scheme: given vector fields
$V_{0},V_{1},V_{2}$ on $\mathbb{R}^{N}$ (notice that H\"{o}rmander's
hypo-ellipticity condition (\ref{formula for greeks}) is not necessarily fulfilled):

\begin{itemize}
\item solve the ordinary differential equations%
\begin{align*}
dY_{s}^{y}(\omega_{1})  &  =\frac{1}{\sqrt{t}}V_{1}(Y_{s}^{y}(\omega
_{1}))ds,\\
dY_{s}^{y}(\omega_{2})  &  =-\frac{1}{\sqrt{t}}V_{1}(Y_{s}^{y}(\omega_{2}))ds,
\end{align*}
which yields flows in the diffusion direction $V_{1}$ starting at
$y\in\mathbb{R}^{N}$.

\item fix the direction for the directional derivative due to the above
construction $e_{1}$ corresponds via formula (\ref{decompo}) to $V_{1}$, hence
$v=\sqrt{t}V_{1}(y)$.

\item the derivative can be approximated through%
\[
\frac{d}{d\epsilon}|_{\epsilon=0}E(f(Y_{t}^{y+\epsilon\sqrt{t}V_{1}%
(y)}))=\frac{1}{2}(f(Y_{t}^{y}(\omega_{1}))-f(Y_{t}^{y}(\omega_{2}%
)))+\mathcal{O}(t^{\frac{3}{2}}),
\]
where the constant in the order estimate depends on the third derivative of
$f$.
\end{itemize}

\end{document}